\documentclass[3p]{elsarticle}
\usepackage{amsmath}
\usepackage{amsfonts,amssymb}

\newtheorem{theorem}{Theorem}{\bf}{\it}
\newtheorem{tab}{Table}{\bf}{\it}
\newtheorem{example}{Example}{\bf}{\rm}
\newtheorem{lemma}{Lemma}{\bf}{\it}
\newtheorem{proposition}{Proposition}{\bf}{\it}

\newcommand{\zr}{\ltimes}
\newcommand{\Real}{\mathbb{R}}
\newcommand{\Co}{\mathbb{C}}

\newcommand{\g}{\mathfrak{g}}
\newcommand{\h}{\mathfrak{h}}

\newcommand{\so}{\mathfrak{so}}
\newcommand{\spin}{\mathfrak{spin}}
\newcommand{\simil}{\mathfrak{sim}}
\def\sl{\mathfrak{sl}}
\newcommand{\gl}{\mathfrak{gl}}
\newcommand{\su}{\mathfrak{su}}
\def\u{\mathfrak{u}}
\def\sp{\mathfrak{sp}}
\newcommand{\f}{\mathfrak{f}}

\newcommand{\R}{\mathcal{R}}
\newcommand{\M}{\mathcal{M}}

\def\P{\mathcal{P}}
\newcommand{\id}{\mathop\text{\rm id}\nolimits}
\newcommand{\spa}{\mathop\text{{\rm span}}\nolimits}
\newcommand{\Hom}{\mathop\text{\rm Hom}\nolimits}
\newcommand{\End}{\mathop\text{\rm End}\nolimits}

\newcommand{\ad}{\mathop\text{\rm ad}\nolimits}
\newcommand{\tr}{\mathop\text{\rm tr}\nolimits}
\newcommand{\pr}{\mathop\text{\rm pr}\nolimits}

\newcommand{\Ric}{\mathop{{\rm Ric}}\nolimits}
\newcommand{\tRic}{\mathop{\widetilde{\rm Ric}}\nolimits}

\begin{document}

\begin{frontmatter}

\title{{One component of the curvature tensor of a Lorentzian manifold}}

\author{Anton S. Galaev}
\address{Department of Mathematics and Statistics, Faculty of Science, Masaryk University in Brno, 
Kotl\'a\v rsk\'a~2, 611~37 Brno,
Czech Republic\\ \ead{galaev@math.muni.cz}}

\begin{abstract}
The holonomy algebra $\g$ of an $n+2$-dimensional Lorentzian manifold $(M,g)$ admitting a parallel distribution of 
isotropic lines is contained in the subalgebra $\simil(n)=(\Real\oplus\so(n))\zr\Real^n\subset\so(1,n+1)$.
An important invariant of $\g$ is its $\so(n)$-projection $\h\subset\so(n)$, which is a Riemannian holonomy algebra.
One component of the curvature tensor of the manifold belongs to the space $\P(\h)$ consisting of linear maps from $\Real^n$ to $\h$ satisfying an identity similar to the Bianchi one. In the present paper the spaces $\P(\h)$ are computed for each possible $\h$. This gives the complete description of the values of the curvature tensor of the manifold $(M,g)$. These results can be applied e.g. to the holonomy classification of the Einstein Lorentzian manifolds.
\end{abstract}


\end{frontmatter}

\section{Introduction} 
The classification of the holonomy algebras of Lorentzian
manifolds is achieved only recently \cite{BB-I,Leistner,Gal5,ESI}.
The most interesting case is when a Lorentzian manifold $(M,g)$ admits a parallel distribution of isotropic lines and the manifold is locally indecomposable, i.e. locally it is not a product of a Lorentzian and a Riemannian manifold. In this case the holonomy algebra $\g$ of $(M,g)$ is contained in
the maximal subalgebra $\simil(n)=(\Real\oplus\so(n))\zr\Real^n$ of the Lorentzian algebra $\so(1,n+1)$ preserving an isotropic line (the dimension of $M$ is $n+2$). There is a number of recent physics literature dealing with these manifolds, see e.g. \cite{BCH,BCH1,CGHP,CFH,CHPP,FF00,G-P,Gibbons09,Hall-Lonie00,LS5}. In particular, in \cite{BCH,BCH1,Gibbons09} expressed the hope that Lorentzian manifolds with the holonomy algebras contained in $\simil(n)$ will found many applications in physics, e.g. they are of interest in M-theory and string theory.

In \cite{Gal1,ESI} the  space $\R(\g)$ of the curvature tensors for each Lorentzian holonomy algebra $\g$, i.e. the space of values of the curvature tensor of a Lorentzian manifold with the holonomy algebra $\g$, are described. Similar results in the Riemannian case \cite{Al} gives a lot of consequences e.g. for Einstein and Ricci-flat manifolds (we explain them below). One component of the space $\R(\g)$ is $$\P(\h)=\{P\in(\Real^n)^*\otimes \h|(P(u)v,w)+(P(v)w,u)+(P(w)u,v)=0\text{ for all }
u,v,w\in\Real^n\},$$ where $\h\subset\so(n)$ is the $\so(n)$-projection of $\g$, which is the holonomy algebra of a Riemannian manifold, and  $(\cdot,\cdot)$ is the inner product on $\Real^n$. In the present paper we compute the spaces $\P(\h)$ for each Riemannian holonomy algebra $\h\subset\so(n)$ (it is enough to assume that $\h\subset\so(n)$ is irreducible).

We introduce the $\h$-equivariant map 
$$\tRic:\P(\h)\to\Real^n, \qquad \tRic(P)=\sum_{i=1}^nP(e_i)e_i,$$ where
$e_1,...,e_n$ is an orthogonal basis of $\Real^n$. If $P$ is a component of the value of the curvature tensor at a point of a manifold $(M,g)$, then $\tRic(P)$ is a component of the Ricci tensor at this point. 
We get the decomposition $$\P(\h)=\P_0(\h)\oplus\P_1(\h),$$ where
$\P_0(\h)$ is the kernel of $\tRic$ and  
$\P_1(\h)$ is its orthogonal complement in $\P(\h)$. If $\h\subset\so(n)$ is irreducible  then
$\P_1(\h)$ is either trivial or it is isomorphic to $\Real^n$. The spaces $\P(\h)$ for
$\h\subset\u(\frac{n}{2})$ are found in \cite{Leistner}. In Section \ref{compP} we compute the spaces $\P(\h)$ for the rest of the algebras. For these computations we turn to the complexification. We consider the representations of semisimple non-simple Lie algebras and the adjoint representations in a unified way. Then we consider case by case the rest of the representations. In particular we show that $\P_0(\h)\neq 0$ and $\P_1(\h)=0$ exactly for the holonomy algebras 
$\su(\frac{n}{2})$, $\sp(\frac{n}{4})$, $\spin(7)$ and $G_2$. Next, $\P_1(\h)\simeq\Real^n$ and $\P_0(\h)\neq 0$ exactly for the holonomy algebras $\so(n)$, $\u(\frac{n}{2})$ and $\sp(\frac{n}{4})\oplus\sp(1)$. For the rest of the Riemannian holonomy algebras (i.e. for the holonomy algebras of symmetric Riemannian spaces different from $\so(n)$, $\u(\frac{n}{2})$ and $\sp(\frac{n}{4})\oplus\sp(1)$) it holds $\P_1(\h)\simeq\Real^n$ and $\P_0(\h)=0$. The result is stated in Table \ref{tabP} of Section \ref{secres}.  We give the explicit forms of some elements $P\in\P(\h)$ in Section \ref{secPespl}.

In \cite{ESI} the study of the holonomy algebras of Einstein Lorentzian manifolds has been begun. It was not possible to complete this study there, as the spaces $\P(\h)$ were not known. In another paper we will use the results obtained here to complete this classification. Here we consider an example dealing with Einstein manifolds admitting a parallel light-like vector field. Any such manifold is Ricci-flat and its holonomy algebra coincides with $\h\zr\Real^n$, where $\h\subset\so(n)$ is the (not necessary irreducible)  holonomy algebra of a Ricci-flat Riemannian manifold.

Necessary facts from the holonomy theory can be found e.g. in \cite{Besse,BCH1,ESI,Jo07}.

Finally remark that the elements of $\P(\h)$ also appear as a component of the curvature tensor of a Riemannian supermanifold $(\M,g)$: it can be checked that $\pr_{\so(T_x\M_{\bar 0})}\circ R_x(\cdot|_{T_x\M_{\bar 0}},\xi)\in\P(\pr_{\so(T_x\M_{\bar 0})}\g)$ for any fixed $\xi\in T_x\M_{\bar 1}$. Here $\g$ is the holonomy algebra of $(\M,g)$ at some point $x$, $T_x\M=T_x\M_{\bar 0}\oplus T_x\M_{\bar 1}$ is the tangent superspace \cite{Galsuphol,GalskBer}. This gives another motivation to the study of the spaces $\P(\h)$.

\section{Preliminaries: the spaces of curvature tensors}\label{curv}
Let $V$ be a vector space and $\g\subset \gl(V)$ a subalgebra. The vector  space $$\R(\g)=\{R\in \Lambda^2
V^*\otimes\g|R(u,v)w+R(v, w)u+R(w, u)v=0 \text{ for all } u,v,w\in V\}$$ is called {\it the space of curvature tensors of
type} $\g$. If there is a pseudo-Euclidean metric $(\cdot,\cdot)$ on $V$ such that $\g\subset\so(V)$, then any $R\in\R(\g)$ satisfies
\begin{equation}\label{symR}  (R(u,v)z,w)= (R(z,w)u,v)\end{equation}
for all $u,v,z,w\in V$.
For $R\in \R(\g)$ define its Ricci tensor by $$\Ric(R)(u,v)=\tr(z\mapsto R(u,z)v),$$ $u,v\in V$.
A subalgebra $\g\subset\gl(V)$ is called {\it a Berger algebra} if $$\g=\spa\{R(u, v)|R\in\R(\g),\,u,v\in
V\},$$ i.e. $\g$ is spanned by the images of the elements $R\in\R(\g)$.
If $(M,g)$ is a pseudo-Riemannian manifold (or, more generally, a manifold $M$ with a torsion-free affine connection $\nabla$) and $\g$ is its holonomy algebra at a point $x\in M$, then $\g$ is a Berger algebra,  the value $R_x$ of the curvature tensor $R$ of $(M,g)$ at the point $x$ belongs to $\R(\g)$, and the value of the Ricci tensor $\Ric$ of $(M,g)$ at the point $x$ coincides with $\Ric(R_x)$. This means that the knowledge 
of the space $\R(\g)$ impose restrictions on the values of $R$ and $\Ric$ and it gives consequences for the geometry of $(M,g)$.
Let us look how does it work with the Einstein condition in the Riemannian case.
The spaces $\R(\h)$ for the holonomy algebras of Riemannian manifolds $\h\subset\so(n)$ are computed by D.~V.~Alekseevsky
in \cite{Al}. Let $\h\subset\so(n)$ be an irreducible Riemannian holonomy algebra.
 The space  $\R(\h)$ admits the following decomposition into
$\h$-modules \begin{equation}\label{razlR}\R(\h)=\R_0(\h)\oplus\R_1(\h)\oplus\R'(\h),\end{equation} where $\R_0(\h)$
consists of the curvature tensors with zero Ricci tensors, $\R_1(\h)$ consists of tensors annihilated by $\h$ (this space is zero or one-dimensional), $\R'(\h)$ is the complement to these two spaces. Each element of $\R'(\h)$ has zero scalar curvature and non-zero Ricci tensor. If $\R(\h)=\R_1(\h)$, then any Riemannian manifold with the holonomy algebra $\h$ is locally symmetric. 
Such subalgebra $\h\subset\so(n)$ is called {\it a symmetric Berger algebra}.  The holonomy algebras of irreducible Riemannian symmetric spaces are exhausted by $\so(n)$, $u(\frac{n}{2})$, $\sp(\frac{n}{4})\oplus\sp(1)$ and by symmetric Berger algebras $\h\subset\so(n)$. Note that $\R(\h)=\R_0(\h)$ if $\h$ is any of the algebras: $\su(\frac{n}{2})$,
$\sp(\frac{n}{4})$, $G_2\subset\so(7)$, $\spin(7)\subset\so(8)$. This implies that each Riemannian manifold with any of
these holonomy algebras is Ricci-flat ($\Ric=0$).
Remark that any locally symmetric Riemannian
manifold is Einstein ($\Ric=\Lambda g$, $\Lambda\in\Real$) and not Ricci-flat.  Next,
$\R(\u(\frac{n}{2}))=\Real\oplus\R'(\u(\frac{n}{2}))\oplus\R(\su(\frac{n}{2}))$ and
$\R(\sp(\frac{n}{4})\oplus\sp(1))=\Real\oplus\R(\sp(\frac{n}{4}))$. Hence any  Riemannian manifold with the holonomy
algebra $\sp(\frac{n}{4})\oplus\sp(1)$ is Einstein and not Ricci-flat, and a Riemannian manifold with the holonomy algebra
$\u(\frac{n}{2})$ can not be Ricci-flat. Finally, if an indecomposable $n$-dimensional Riemannian manifold is Ricci-flat,
then its holonomy algebra is one of $\so(n)$, $\su(\frac{n}{2})$, $\sp(\frac{n}{4})$,  $G_2\subset\so(7)$,
$\spin(7)\subset\so(8)$. 

Irreducible holonomy algebras $\g\subset\gl(n,\Real)$ of torsion-free affine connections are classified by S.~Merkulov and L.~Schwachh\"ofer in \cite{M-Sch,Sch}.
In \cite{Armstrong} S.~Armstrong, analysing the spaces $\R(\g)$, found which of these holonomy algebras correspond to Ricci-flat connections.

Consider now the case of Lorentzian manifolds.  From the Wu  Theorem \cite{Wu} it follows that any Lorentzian manifold $(M,g)$ is either locally a product of the manifold $(\Real,-(dt)^2)$ and of a Riemannian manifold, or of a Lorentzian and a Riemannian manifold, or it is locally indecomposable,
i.e. it does not admit such decompositions. If the manifold $(M,g)$ is simply
connected and geodesically complete, then these decompositions are global. This allows us to consider locally
indecomposable Lorentzian manifolds. The only irreducible Lorentzian holonomy algebra is the whole Lie algebra
$\so(1,n+1)$ \cite{Ber} (the dimension of $M$ is $n+2$). If $(M,g)$ is locally indecomposable and its holonomy  algebra $\g$ is different from $\so(1,n+1)$, then $\g$ preserves an isotropic line of the tangent space and $(M,g)$ locally admits   parallel distributions of isotropic lines.
If $M$ is simply connected, then there exists a global  parallel distribution of isotropic lines.

\begin{example}\label{ex1} Let $(M,g)$ be an $n+2$-dimensional locally indecomposable Lorentzian manifold admitting a parallel light-like vector field $X$.  Let $x\in M$. We identify 
the tangent space $T_xM$ with the Minkowski space $\Real^{1,n+1}$. Let $p\in\Real^{1,n+1}$ be the value of $X$ at the point $x$. Choose a basis $p,e_1,...,e_n,q$ of $\Real^{1,n+1}$ with the following non-zero values of $g_x$: $g_x(p,q)=1$, $g_x(e_i,e_i)=1$. The subalgebra of $\so(1,n+1)$ preserving $p$
has the form $\so(n)\zr(p\wedge \Real^n)$ (we identify   $\Real^n$ with $\spa\{e_1,...,e_n\}$, and  $\wedge^2\Real^{1,n+1}$ with $\so(1,n+1)$ such that it holds $(u\wedge v)(z)=(u,z)v-(v,z)u$; in particular, $\wedge^2\Real^{n}$ is identified with $\so(n)$).
Then the holonomy algebra of $(M,g)$ at the point $x$ is contained in the algebra $\h\zr(p\wedge \Real^n)$, where $\h\subset\so(n)$ is the 
$\so(n)$-projection of $\g$, which is a (not necessary irreducible) Riemannian holonomy algebra \cite{Leistner}. The value $R_x$ satisfies 
$$ R_x(p,\cdot)=0,\quad R_x(u,v)=R_0(u,v)+p\wedge (P(u)v-P(v)u),\quad R_x(u,q)=P(u)-p\wedge T(u)$$ for all $u,v\in \Real^n$.
Here $R_0\in\R(\h)$, $P\in\P(\h)$, and $T\in\End(\Real^n)$, $T^*=T$ \cite{Gal1,ESI}.
In particular, if $\g=\h\zr(p\wedge \Real^n)$, then $\R(\g)\simeq\R(\h)\oplus\P(\h)\oplus\odot^2\Real^n$. Thus the only unknown space in this decomposition is $\P(\h)$. The spaces $\R(\g)$ for other Lorentzian holonomy algebras have similar description \cite{Gal1,ESI}. 
\end{example}

\section{Main result}\label{secres}

Now we begin to study the space $\P(\h)$, where $\h\subset\so(n)$ is an irreducible subalgebra. Consider the $\h$-equivariant map
$$\tRic:\P(\h)\to\Real^n, \qquad \tRic(P)=\sum_{i=1}^nP(e_i)e_i.$$ This definition does not depend on the choice of the
orthogonal basis $e_1,...,e_n$ of $\Real^n$. Denote by $\P_0(\h)$ the kernel of $\tRic$ and let $\P_1(\h)$ be its
orthogonal complement in $\P(\h)$. Thus, $$\P(\h)=\P_0(\h)\oplus\P_1(\h).$$ Since $\h\subset\so(n)$ is irreducible and the
map $\tRic$ is $\h$-equivariant, $\P_1(\h)$ is either trivial or isomorphic to $\Real^n$. The spaces $\P(\h)$ for
$\h\subset\u(\frac{n}{2})$ are found in \cite{Leistner}. In Section \ref{compP} we compute the spaces $\P(\h)$ for the
remaining Riemannian holonomy algebras. The result is given in Table \ref{tabP} (for a compact Lie algebra $\h$,
$V_\Lambda$ denotes the irreducible representation of $\h$ given by the irreducible representation of $\h\otimes
\mathbb{C}$ with the highest weight $\Lambda$; $((\odot^2(\Co^m)^*\otimes\Co^m)_0$ denotes the subspace of $\odot^2(\Co^m)^*\otimes\Co^m$ consisting of tensors such that the contraction of the upper index with any down index gives zero). 

\begin{tab}\label{tabP} The spaces $\P(\h)$ for irreducible Riemannian holonomy algebras $\h\subset\so(n)$.
$$\begin{array}{|c|c|c|c|} \hline \h\subset\so(n)&\P_1(\h)&\P_0(\h)&\dim\P_0(\h)\\\hline \so(2)&\Real^2&0&0\\\hline \so(3)&\Real^3&V_{4\pi_1}&5\\\hline
\so(4)&\Real^4&V_{3\pi_1+\pi_1'}\oplus V_{\pi_1+3\pi_1'}&16\\\hline \so(n),\,n\geq 5&\Real^n&V_{\pi_1+\pi_2}&\frac{(n-2)n(n+2)}{3}\\\hline \u(m),\,
n=2m\geq 4&\Real^n&(\odot^2(\Co^m)^*\otimes\Co^m)_0&m^2(m-1)\\\hline \su(m),\, n=2m\geq 4&0&(\odot^2(\Co^m)^*\otimes\Co^m)_0&m^2(m-1)\\\hline
\sp(m)\oplus\sp(1),\, n=4m\geq 8&\Real^n&\odot^3(\Co^{2m})^*&\frac{m(m+1)(m+2)}{3}\\\hline \sp(m),\, n=4m\geq 8&0&\odot^3(\Co^{2m})^*&\frac{m(m+1)(m+2)}{3}\\\hline
G_2\subset\so(7)&0&V_{\pi_1+\pi_2}&64\\\hline \spin(7)\subset\so(8)&0&V_{\pi_2+\pi_3}&112\\\hline \h\subset\so(n),\,n\geq 4,&\Real^n&0&0\\ \text{ is a symmetric Berger algebra}&&&\\\hline\end{array}$$\end{tab}

\begin{example} In the settings of Example \ref{ex1} for the value $\Ric_x$ of the Ricci tensor we have 
\begin{align*}\Ric_x(p,\cdot)&=0, \quad \Ric_x(u,v)=\Ric(R_0)(u,v),\\
 \Ric_x(u,q)&=g_x(u,\tRic(P)),\quad \Ric_x(q,q)=\tr T.\end{align*} 
Suppose that $(M,g)$ is an Einstein  manifold, i.e. $g$ satisfies the equation of General Relativity in the absence of  
matter $$\Ric=\Lambda g,\quad\Lambda\in\Real.$$
 Using the expression for $\Ric_x$ in \cite{ESI} it is proved
that $(M,g)$ is Ricci-flat (i.e. $\Lambda=0$)
 and $\g=\h\zr(p\wedge \Real^n)$. Using the above results for the space $\P(\h)$
it can be shown that $\h\subset\so(n)$ is the holonomy algebra of a Ricci-flat Riemannian manifold, in particular
there are decompositions \begin{equation}\label{razlAB}\Real^{n}=\Real^{n_1}\oplus\cdots\oplus\Real^{n_s}\oplus\Real^{n_{s+1}},\quad \h=\h_1\oplus\cdots\oplus\h_s\oplus\{0\}\end{equation} such that $\h_i(\Real^{n_j})=0$ for $i\neq j$, and each 
$\h_i\subset\so(n_i)$ coincides with one of the Lie algebras
$\so(n_i)$, $\su(\frac{n_i}{2})$, $\sp(\frac{n_i}{4})$, $G_2\subset\so(7)$, $\spin(7)\subset\so(8)$. 
\end{example}

Note that any Riemannian holonomy algebra $\h\subset\so(n)$ admits the decomposition \eqref{razlAB} and it holds $$\R(\h)=\R(\h_1)\oplus\cdots\oplus\R(\h_s),\qquad \P(\h)=\P(\h_1)\oplus\cdots\oplus\P(\h_s).$$

Consider the natural $\h$-equivariant map $$\tau:\Real^n\otimes\R(\h)\to\P(\h),\qquad \tau(u\otimes R)=R(\cdot,u).$$

\begin{theorem} For any irreducible subalgebra $\h\subset\so(n)$, $n\geq 4$,  the $\h$-equivariant map
$\tau:\Real^n\otimes\R(\h)\to\P(\h)$ is surjective. Moreover, $\tau(\Real^n\otimes\R_0(\h))=\P_0(\h)$ and
$\tau(\Real^n\otimes\R_1(\h))=\P_1(\h)$.\end{theorem}

{\bf Proof.} If $\h\subset\so(n)$ is a Riemannian holonomy algebra, then the theorem  follows from Table \ref{tabP} and
the results from \cite{Al}. We claim that the only irreducible subalgebra $\h\subset\so(n)$ with $\R(\h)\neq 0$ that is
not a Berger algebra is $\sp(\frac{n}{4})\oplus\Real J$, where $J$ is a complex structure on $\Real^n$. Similarly, the
only irreducible subalgebra $\h\subset\so(n)$ with $\P(\h)\neq 0$ that is not spanned by images of the elements from $\P(\h)$ is
$\sp(\frac{n}{4})\oplus\Real J$. It is enough to prove the second claim. Let $\h\subset\so(n)$ be an irreducible
subalgebra  with $\P(\h)\neq 0$ and such that $\h\subset\so(n)$ is not spanned by images of the elements from $\P(\h)$. If
$\h\not\subset\u(\frac{n}{2})$, then $\h\otimes\Co$ is semisimple non-simple Lie algebra, and
$\h\otimes\Co\subset\so(n,\Co)$ is an irreducible subalgebra with $\P(\h\otimes\Co)\neq 0$. Let $\h_1\subset\h\otimes\Co$
be the ideal generated by the images of the elements from $\P(\h\otimes\Co)$. Proposition \ref{prop2} below shows that
$\h\otimes\Co\subsetneq\sl(2,\Co)\oplus \sp(\frac{n}{2},\Co)$ and $\sl(2,\Co)\subset\h\otimes\Co$. The proof of
Proposition \ref{prop2} implies that $\sl(2,\Co)\subset \h_1$. We get that there is a proper ideal
$\h_2\subset\h\otimes\Co$ such that $\h\otimes\Co=\h_1\oplus\h_2$ and $\P(\h\otimes\Co)=\P(\h_1)$. The proof of
Proposition \ref{prop2} shows that this is impossible. Hence, $\h\subset\u(\frac{n}{2})$. From the results of
\cite{Leistner} it follows that $\h=\sp(\frac{n}{4})\oplus\Real J$. Finally note that $\R(\sp(\frac{n}{4})\oplus\Real
J)=\R(\sp(\frac{n}{4}))$ and $\P(\sp(\frac{n}{4})\oplus\Real J)=\P(\sp(\frac{n}{4}))$.  \qed

\section{The explicit form of some $P\in\P(\h)$}\label{secPespl}
Using the above results and results from \cite{Al}, we can now explicitly give the spaces $\P(\h)$ in some cases.

From the results of \cite{Leistner} it follows that $\P(\u(m))\simeq\odot^2(\Co^m)^*\otimes\Co^m$. Let us describe this
isomorphism  in the following way. Let $S\in \odot^2(\Co^m)^*\otimes\Co^m\subset(\Co^m)^*\otimes \gl(m,\Co)$. We fix an
identification $\Co^m=\Real^{2m}=\Real^m\oplus i\Real^m$ and choose a basis $e_1,...,e_m$ of $\Real^m$. Define the complex
numbers $S_{abc}$, $a,b,c=1,...,m$ such that $S(e_a)e_b=\sum_c S_{acb}e_c$. It holds $S_{abc}=S_{cba}$. Define a map
$S_1:\Real^{2m}\to\gl(2m,\Real)$ by the conditions $S_1(e_a)e_b=\sum_c \overline {S_{abc}}e_c$, $S_1(ie_a)=-iS_1(e_a)$,
and $S_1(e_a)ie_b=iS_1(e_a)e_b$. It is easy to check that $P=S-S_1:\Real^{2m}\to\gl(2m,\Real)$ belongs to $\P(\u(n))$ and
any element of $\P(\u(n))$ is of this form. Such element belongs to $\P(\su(n))$ if and only if $\sum_b S_{abb}=0$ for all
$a=1,...,m$, i.e. $S\in (\odot^2(\Co^m)^*\otimes\Co^m)_0$. If $m=2k$, i.e. $n=4k$,  then $P$ belongs to $\P(\sp(k))$ if
and only if $S(e_a)\in\sp(2k,\Co)$, $a=1,...,m$, i.e. $S\in(\sp(2k,\Co))^{(1)}\simeq\odot^3(\Co^{2k})^*$.

In \cite{ESI} it is shown that any $P\in\P(\u(m))$ satisfies $ (\tRic(P),x)=-\tr_\Co P(Jx)$ for all $x\in\Real^{2m}$.

Note that $\R_1(\so(n))\oplus\R'(\so(n))\simeq \odot^2\Real^n$.
Any $R\in \R_1(\so(n))\oplus\R'(\so(n))$ is of the
form $R=R_S$, where $S:\Real^{n}\to\Real^{n}$ is a symmetric linear map  and $$R_S(x,y)=Sx\wedge y+
x\wedge Sy.$$ Similarly, $\R_1(\so(n))$ is spanned by the element $R=R_{\frac{\id}{2}}$, i.e. $R(x,y)=x\wedge y$.
Next, $\tau(\Real^n,\R_1(\so(n))\oplus\R'(\so(n)))=\P(\so(n))$. Hence $\P(\so(n))$ is spanned by the
elements $P$ of the form $$P(y)=Sy\wedge x+y\wedge Sx,$$ where $x\in\Real^{n}$ and $S\in\odot^2\Real^n$  are  fixed, and
$y\in\Real^{n}$ is any vector. For such $P$ it holds $\tRic(P)=(\tr S-S)x$. We conclude that $\P_0(\so(n))$ is spanned by
the elements $P$ of the form $$P(y)=Sy\wedge x,$$ where $x\in\Real^{n}$ and $S\in\odot^2\Real^n$ satisfy $\tr S=0$ and
$Sx=0$, and $y\in\Real^{n}$ is any vector.

The isomorphism $\P_1(\so(n))\simeq\Real^n$ is defined in the following way: $x\in\Real^n$ corresponds to
$P=x\wedge\cdot\in\P_1(\so(n))$, i.e. $P(y)=x\wedge y$ for all $y\in\Real^n$.

Any $P\in\P_1(\u(m))$ has the form $$P(y)=-\frac{1}{2} (Jx,y)J+\frac{1}{4}(x\wedge y+Jx\wedge Jy),$$  where $J$ is the
complex structure on $\Real^{2m}$, $x\in\Real^{2m}$ is fixed, and $y\in\Real^{2m}$ is any vector.

Any $P\in\P_1(\sp(m)\oplus\sp(1))$ has the form $$P(y)= -\frac{1}{2}\sum_{\alpha=1}^3g(J_\alpha
x,y)J_\alpha+\frac{1}{4}\big(x\wedge y+ \sum_{\alpha=1}^3J_\alpha x\wedge J_\alpha y\big),$$ where $(J_1,J_2,J_3)$ is the
quaternionic structure on $\Real^{4m}$, $x\in\Real^{4m}$ is fixed, and $y\in\Real^{4m}$ is any vector.

We will see that for $\h\subset\so(\h)$, where $\h$ is a compact simple Lie algebra any $P\in\P(\h)=\P_1(\h)$ has the form
$$P(y)=[x,y],$$ where $x\in\h$ is fixed, and $y\in\h$ is any element. 

If $\h\subset\so(n)$ is a symmetric Berger algebra, then $\P(\h)=\{R(\cdot,x)|\,x\in\Real^n\}$, where $R$ is a generator of $\R(\h)\simeq\Real$. 

In general, let $\h\subset\so(n)$ be an irreducible subalgebra and $P\in\P_1(\h)$. Then
$\tRic(P)\wedge\cdot\in\P_1(\so(n))$. Furthermore, it is easy to check that
$\tRic\big(P+\frac{1}{n-1}\tRic(P)\wedge\cdot\big)=0$, that is $P+\frac{1}{n-1}\tRic(P)\wedge\cdot\in\P_0(\so(n))$. Thus
the inclusion $\P_1(\h)\subset\P(\so(n))=\P_0(\so(n))\oplus\P_1(\so(n))$ is given by $$P\in
P_1(\h)\mapsto\big(P+\frac{1}{n-1}\tRic(P)\wedge\cdot,-\frac{1}{n-1}\tRic(P)\wedge\cdot\big)\in\P_0(\so(n))\oplus\P_1(\so(n)).$$
This construction defines  the tensor $W=P+\frac{1}{n-1}\tRic(P)\wedge\cdot$, which is the analog of the Weyl tensor for
$P\in\P(\h)$, and this tensor is a component of the Weyl tensor of a Lorentzian manifold.

\section{Computation of the spaces $\P(\h)$}\label{compP}

Let $\h\subset\so(n)$ be an irreducible Riemannian holonomy algebra. Since for the subalgebras $\h\subset\u(\frac{n}{2})$
the spaces $\P(\h)$ are found in \cite{Leistner}, we may assume that $\h\not\subset\u(\frac{n}{2})$, then the subalgebra
$\h\otimes\Co\subset\so(n,\Co)$ is irreducible and it is enough to find the space $\P(\h\otimes\Co)$, which equals
$\P(\h)\otimes\Co$.


Here we compute the spaces $\P(\h)$ for all irreducible Berger subalgebras $\h\subset\so(n,\Co)$. The only non-symmetric
Berger subalgebras  $\h\subset\so(n,\Co)$ (i.e. subalgebras with $\R(\h)\neq\R_1(\h)$) are $\so(n,\Co)$,
$\sl(2,\Co)\oplus\sp(2m,\Co)\subset\so(4m,\Co)$ ($m\geq 2$), $\spin(7,\Co)\subset \so(8,\Co)$ and
$G^\Co_2\subset\so(7,\Co)$. The symmetric Berger subalgebras $\h\subset\so(m,\Co)$ (i.e. subalgebras with
$\R(\h)=\R_1(\h)$) are given in Table \ref{symBer} taken from \cite{Sch}.

For simple Lie algebras we use the notation from \cite{V-O}. For some computations we use the package LiE \cite{LiE}.
Remark  that the numbering of the vertices on the Dynkin diagrams for some simple Lie algebras in \cite{V-O} and
\cite{LiE} are different.

\begin{tab} \label{symBer} Irreducible symmetric Berger subalgebras $\h\subset\so(m,\Co)=\so(V)$.
$$\begin{array}{|c|c|c|}\hline No.&\h&V\\\hline 1&\sp(2n,\Co),\,n\geq 3&V_{\pi_2}=\Lambda^2\Co^{2n}/\Co\omega\\\hline
2&\so(n,\Co),\,n\geq 3,\,n\neq 4&V_{2\pi_1}=\odot^2\Co^{n}/\Co g \\\hline 3&\h \text{ is a simple Lie algebra} &\h\\\hline
4&\so(9,\Co)& (\Delta_9)^\Co\\\hline 5&\sp(8,\Co)&V_{\pi_4}=\Lambda^4\Co^{8}/(\omega\wedge\Lambda^2\Co^{8})\\\hline
6&F_4^\Co & V_{\pi_1}=\Co^{26}\\\hline 7&\sl(8,\Co)&V_{\pi_4}=\Lambda^4\Co^{8}\\\hline 8&\so(16,\Co)&
(\Delta_{16}^+)^\Co\\\hline 9&\so(p,\Co)\oplus \so(q,\Co),\, p,q\geq 3&\Co^{p}\otimes\Co^{q}\\\hline 10&\sp(2p,\Co)\oplus
\sp(2q,\Co),\,p,q\geq 2&\Co^{2p}\otimes\Co^{2q}\\\hline 11&\sl(2,\Co)\oplus
\sl(2,\Co)&\Co^{2}\otimes\odot^3\Co^{2}\\\hline
12&\sp(6,\Co)\oplus \sl(2,\Co)&V_{\pi_3}\otimes\Co^2=(\Lambda^3\Co^{6}/(\omega\wedge\Co^{6}))\otimes\Co^2\\\hline
13&\sl(6,\Co)\oplus \sl(2,\Co)&V_{\pi_3}\otimes\Co^2=\Lambda^3\Co^{6}\otimes\Co^2\\\hline 14&\so(12,\Co)\oplus \sl(2,\Co)&
(\Delta_{12}^+)^\Co\otimes\Co^2\\\hline 15&E_7^\Co\oplus \sl(2,\Co)&V_{\pi_1}\otimes\Co^2= \Co^{56}\otimes\Co^2\\\hline
\end{array} $$
\end{tab}

\begin{lemma}\label{lemRP} Let $V$ be a real or complex vector space and $\h\subset\so(V)$. Then a linear map $R:\Lambda^2
V\to\h$ belongs to $\R(\h)$ if and only if for each $x\in V$ it holds $R(\cdot,x)\in\P(\h)$.\end{lemma}

{\it Proof.} If $R\in\R(\h)$, then the inclusion $R(\cdot,x)\in\P(\h)$ follows from \eqref{symR} and the Bianchi identity.
Conversely, if $R(\cdot,x)\in\P(\h)$ for each $x\in V$, then it is not hard to prove that $R$ satisfies \eqref{symR}, and
using this it is easy to see that $R\in\R(\h)$.  \qed

Remark that the above lemma can be also applied for irreducible submodules $U\subset\R(\h)$.

\begin{lemma}\label{Slem1} Let $\h\subset\so(n,\Co)=\so(V)$ be an irreducible subalgebra. Then the decomposition of the tensor product $V\otimes \h$ into irreducible $\h$-modules is of the form
$V\otimes \h=kV\oplus (\oplus_\lambda V_\lambda)$, where  $k$ is the number of non-zero labels on the Dynkin diagram for
the representation of $\h$ on $V$, and $V_\lambda$ are pairwise non-isomorphic irreducible $\h$-modules that are not
isomorphic to $V$. \end{lemma}

{\it Proof.} The number of irreducible submodules isomorphic to the highest weight module $V_\lambda$ in the product
$V\otimes \h$ is equal to $$\dim\{v\in\h_{\lambda-\Lambda}|(\ad_{A_i})^{\Lambda_i+1}v=0,\,\, i=1,...,l\},$$ where
$\Lambda$ is the highest weight of $V$, $l$ is the rang of $\h$, $A_i$ are canonical generators of $\h$ corresponding to
the simple positive roots, and $\Lambda_i$ are the labels on the Dynkin diagram defining $\Lambda$, see e.g. \cite{V-O}.
If $\Lambda\neq\lambda$, then $\dim\h_{\lambda-\Lambda}$ equals either $0$ or $1$. This shows that all $V_\lambda$ are
pairwise different. We get that $$k=\dim\{v\in\h_{0}|(\ad_{A_i})^{\Lambda_i+1}v=0,\,\, i=1,...,l\},$$ where $\h_{0}$ is
the fixed Cartan subalgebra of $\h$. If $\Lambda_i>0$, then obviously $(\ad_{A_i})^{\Lambda_i+1}v=0$. We get
$$k=\dim\{v\in\h_{0}|[A_i,v]=0\text{ whenever } \Lambda_i=0\}.$$ The matrix of the obtained homogeneous system of linear
equations consists of the lines of the Cartan matrix of $\h$ corresponding to $i$ with $\Lambda_i=0$. Since the Cartan
matrix is non-degenerate, we immediately get the proof of the lemma.  \qed

\vskip0.5cm

\begin{proposition}\label{prop1} Let $V$ be a real or complex vector space and $\h\subset\so(V)$ be an irreducible subalgebra.
If $V\otimes\h$ contains only one irreducible submodule isomorphic to $V$, then $\P_1(\h)\simeq V$ if and only if
$\R_1(\h)\simeq\mathbb{F}$, where $\mathbb{F}=\Real$ or $\Co$, respectively.
\end{proposition}

{\bf Proof.} Let $\Hom_0(V,\h)\subset\Hom(V,\h)$ be the subset consisting of the maps $\varphi:V\to\h$ such that
$\sum_{i=1}^n\varphi(e_i)e_i=0$. Denote by $\Hom_1(V,\h)$ its orthogonal complement, then $$\Hom(V,\h)=\Hom_0(V,\h)\oplus
\Hom_1(V,\h).$$ It is easy to see that $\Hom_1(V,\h)\simeq V$. Note that $\Hom_1(V,\so(n))=\P_1(\so(n))$.

\begin{lemma} $\Hom_1(V,\h)=\{\pr_\h\circ P|P\in\P_1(\so(n))\}$.\end{lemma}
Let $\varphi\in\Hom_0(V,\h)$ and $P\in\P_1(\so(n))$, then $$(\varphi,\pr_\h\circ P)=\sum_{i,j=1}^n(e_i\otimes
\varphi(e_i),e_j\otimes\pr_\h\circ P(e_j))=\sum_{i=1}^n(\varphi(e_i),\pr_\h\circ
P(e_i))=\sum_{i=1}^n(\varphi(e_i),P(e_i))=(\varphi,P)=0,$$ where we used the scalar products on different tensor
spaces and the fact that $\Hom_0(V,\h)\subset \Hom_0(V,\so(n))$ is orthogonal to $\Hom_1(V,\so(n))$. On the other hand,
suppose that $\pr_\h\circ P=0$. Recall that $P$ is of the form $x_0\wedge \cdot$ for some $x_0\in V$. Then $\h$
annihilates $x_0$ and we get $x_0=0$. The lemma is proved.  \qed

Note that $\odot^2\h$ contains $\id_\h$. Moreover, either $\R_1(\h)=0$, or $\R_1(\h)=\mathbb{F} \id_\h$. Let $R:\Lambda^2
V\to\h$ be the extension of $\id_\h$ such that $R|_{\h^\perp}=0$, then $R(x,y)=\pr_\h(x\wedge y)$ for all $x,y\in V$. It
is clear that either $\P_1(\h)=0$, or $\P_1(\h)=\Hom_1(V,\h)$.  Lemma \ref{lemRP} implies that $\P_1(\h)=\Hom_1(V,\h)$ if
and only if $\R_1(\h)=\mathbb{F} \id_\h$.  \qed

\begin{proposition}\label{prop2} Let $\h_1\subset\gl(V_1)$ and $\h_2\subset\gl(V_2)$ be irreducible complex subalgebras and
$\h=\h_1\oplus\h_2\subset\so(V_1\otimes V_2)=\so(V)$. If   $\h$ is different from $\sl(2,\Co)\oplus\sp(2m,\Co)$, then
$\P_0(\h)=0$. Consequently, if $\h$ is a symmetric Berger algebra, then $\P(\h)=\P_1(\h)\simeq V$. Moreover,
$\P_1(\sl(2,\Co)\oplus\sp(2m,\Co))\simeq V$ and $\P_0(\sl(2,\Co)\oplus\sp(2m,\Co))=(\sp(2m,\Co))^{(1)}\oplus
(\sp(2m,\Co))^{(1)}$, where $(\sp(2m,\Co))^{(1)}\simeq \odot^3(\Co^{2m})^*$ is the first prolongation of the subalgebra
$\sp(2m,\Co)\subset\gl(2m,\Co)$.\end{proposition}

{\it Proof.} First suppose that the dimensions of $V_1$ and $V_2$ are greater then 2. Let $\h$ be one of the following:
$\h=\so(V_1)\oplus\so(V_2)$,  $\h=\sp(V_1)\oplus\sp(V_2)$.

We claim that there is no $P\in\P(\h)$ taking values either in $\h_1$ or in $\h_2$, i.e.
$\P(\h_1\subset\so(V))=\P(\h_2\subset\so(V))=0$. Indeed, since the dimensions of $V_1$ and $V_2$ are greater then 2,  both
$\h_1$ and $\h_2$ preserve more then two vector subspaces of $V$ and act in these subspaces in the same time. From this it
is easily follows that $\P(\h_1\subset\so(V))=\P(\h_2\subset\so(V))=0$. The claim is proved.

Note that $\h\subset\so(V)$ is a symmetric Berger  subalgebra, and  there are exactly two non-zero labels on the Dynkin
diagram of $\h$ defining the representation $V$. One of these labels is on the Dynkin diagram of $\h_1$ and the other one
is on the Dynkin diagram of $\h_2$. Hence $V\otimes \h$ contains two irreducible components isomorphic to $V$. Next,
$V\otimes \h= (\h_1\otimes V_1\otimes V_2)\oplus (\h_2\otimes V_1\otimes  V_2)$. This shows that one irreducible component
$V\subset V\otimes\h$ belongs two $\h_1\otimes V_1\otimes V_2$ and another one belongs to $\h_2\otimes V_1\otimes  V_2$.
Hence none of them belong two $\P(\h)$. On the other side, $\P(\h)$ contains $\P_1(\h)\simeq V$. We conclude that
$(V\oplus V)\cap \P(\h)=\P_1(\h)\simeq V$. If $V_\lambda\subset V\otimes\h$ is an irreducible $\h$-submodule not
isomorphic to $V$, then $V_\lambda$ is contained ether in $\h_1\otimes V_1\otimes V_2$, or in $\h_2\otimes V_1\otimes
V_2$, i.e. it is not contained in $\P(\h)$.

If $\h=\h_1\oplus\h_2\subset\so(V_1\otimes V_2)=\so(V)$ is an irreducible subalgebra different from the above two, then it
is properly contained either in $\f=\so(V_1)\oplus\so(V_2)$, or in $\f=\sp(V_1)\oplus\sp(V_2)$. Next,
 $\P(\h)=(V\otimes\h)\cap\P(\f)$. Since $\P(\f)\simeq V$ is an irreducible $\h$-module and the images of the elements of $\P(\f)$ span $\f$, we get that $\P(\h)=0$.

If $\dim V_1=2$, then $\h_1=\sl(2,\Co)$ and $\h_2\subset\sp(2m,\Co)$ is a proper irreducible subalgebra, $m\geq 2$.
Consider the Lie algebra $\h=\sl(2,\Co)\oplus\sp(2m,\Co)$. It is easy to see that $\P(\sp(2m,\Co)\subset\so(V))
=(\sp(2m,\Co))^{(1)}\oplus(\sp(2m,\Co))^{(1)}$.  Using this and the above arguments, we get that $\P_1(\h)\simeq V$ and
$\P_0(\h)=(\sp(2m,\Co))^{(1)}\oplus (\sp(2m,\Co))^{(1)}$. If $\h_2\subset\sp(2m,\Co)$ is a proper subalgebra, then
$\P_0(\sl(2,\Co)\oplus \h_2)=(V\otimes (\sl(2,\Co)\oplus \h_2))\cap \P_0(\sl(2,\Co)\oplus \sp(2m,\Co))$ and this
intersection is zero, since it holds $(\h_2)^{(1)}=0$ \cite{Sch}.   The proposition is proved.  \qed

\vskip0.5cm

Let $\h$ be simple and $\delta$ be its highest root. Let $V=V_\Lambda$. The tensor product $V\otimes \h$ contains the
$\h$-submodule $V_{\Lambda+\delta}$. Let $A_\delta\in\h$ and $v_\Lambda\in V$ be the highest root and highest weight
vectors, respectively. Then $v_\Lambda\otimes A_\delta\in V_{\Lambda+\delta}\subset V\otimes\h$, and
$V_{\Lambda+\delta}\subset \P(\h)$ if and only if $v_\Lambda\otimes A_\delta\in \P(\h)$. It is easy to see that the
condition $v_\Lambda\otimes A_\delta\in \P(\h)$ holds if and only if $A_\delta$ has rank two and $A_\delta
v_{-\Lambda}\neq 0$, where $v_{-\Lambda}$ is the lowest vector in $V$ (note that $v_{\Lambda}$ and $v_{-\Lambda}$ are
isotropic, and $ (v_{\Lambda},v_{-\Lambda})\neq 0$ ). If $A_\delta$ has rank two, then $A_\delta\odot
A_\delta\in\R(\h)$ and $\h\subset\so(V)$ is a non-symmetric Berger subalgebra. We have proved that if $\h\subset\so(V)$ is
a symmetric Berger subalgebra or $\R(\h)=0$, then $V_{\Lambda+\delta}\cap \P(\h)=0$. Thus we need only to consider
irreducible submodules $V_\lambda\subset V\otimes\h$ not isomorphic to  $V_{\Lambda+\delta}$ and $V$ (if the representation is given by the Dynkin diagram with only one non-zero label). We will see that such
submodules are never contained in $\P(\h)$.

\vskip0.5cm

It is easy to get that $\Co^n\otimes\so(n,\Co)=\Co^n\oplus V_{\pi_1\oplus\pi_2}\oplus V_{\pi_3}$ ($n\geq 5$),
$\Co^7\otimes G^\Co_2=\Co^7\oplus V_{\pi_1+\pi_2}\oplus V_{2\pi_2}$ and $\Co^8\otimes \so(7)=\Co^8\oplus
V_{\pi_2+\pi_3}\oplus V_{\pi_1+\pi_3}$.
\begin{lemma} We have $\P_1(\so(n,\Co))\simeq\Co^n$, $\P_0(\so(n,\Co))= V_{\pi_1\oplus\pi_2}$ ($n\geq 5$),
$\P(G^\Co_2)=\P_0(G^{\Co}_2)= V_{\pi_1+\pi_2}$ and $\P(\spin(7,\Co))=\P_0(\spin(7,\Co))= V_{\pi_2+\pi_3}$.\end{lemma}

{\it Proof.} From Proposition \ref{prop1} it follows that $\P_1(\so(n,\Co))\simeq\Co^n$, and for both Lie algebras
$G_2^\Co\subset \so(7,\Co)$ and $\spin(7,\Co)\subset\so(8,\Co)$ it holds $\P_1(\h)=0$, i.e. $\P(\h)\cap V=0$.

Let $\h=G^\Co_2$. We have $\Lambda=\pi_1=\epsilon_1$ and $\delta=\pi_2=\epsilon_1-\epsilon_3$. It is easy to see that
$v_\Lambda\otimes A_\delta\in \P(\h)$, i.e. $V_{\Lambda+\delta}\subset \P(\h)$. Next, $\subset V\otimes \h$ contains a
3-dimensional vector subspace of weight $2\pi_2$. This subspace is spanned by the vectors $v_{\epsilon_1}\otimes
A_{\epsilon_1}$, $v_{-\epsilon_2}\otimes A_{\epsilon_1-\epsilon_3}$ and $v_{-\epsilon_3}\otimes
A_{\epsilon_1-\epsilon_2}$, where $A_\mu$ denotes a non-zero root element in $\h$ of weight $\mu$. Moreover this subspace
has a 1-dimensional intersection with $V_{2\pi_2}$ and a 2-dimensional intersection with $V_{\pi_1+\pi_2}$. This shows
that $V_{2\pi_2}\subset \P(\h)$ if and only if all these three vectors belong to $\P(\h)$. To see that
$v_{\epsilon_1}\otimes A_{\epsilon_1}\not\in \P(\h)$ it is enough to write down the definition of $P\in\P(\h)$ for the vectors
$v_{-\epsilon_1}$, $v_{\epsilon_2}$ and $v_{\epsilon_3}$.

The Lie algebras $\so(n,\Co)$ and $\spin(7,\Co)\subset\so(8,\Co)$  can be considered in the same way.  \qed

\begin{lemma} Let $\h$ be a simple complex Lie algebra different from $\sl(2,\Co)$. Then for the adjoint representation $\h\subset\so(\h)$ we have $\P(\h)=\P_1(\h)\simeq \h$. \end{lemma}
{\it Proof.} Since the space $\R(\h)$ is one-dimensional and it is spanned by the Lie brackets of $\h$ \cite{Sch}, any
element $P\in\P_1(\h)$ is of the form $P(\cdot)=[\cdot,x]$ for some $x\in\h$.

Denote by $(\cdot,\cdot)$ the Killing form on $\h$.  Let $P\in\P(\h)$, then for any $x,y,z\in\h$ it holds
 \begin{multline*} ([P(x),y],z)+ ([x,P(y)],z)=- ([P(y),z],x)- ([P(z),x],y)+ ([x,P(y)],z)\\
 =- ([P(z),x],y)=- (P(z),[x,y]).\end{multline*}
Since $\h\otimes\h=\odot^2\h\oplus\Lambda^2\h$ and the adjoint representation of any simple $\h$ different from
$\sl(n,\Co)$ is given by the Dynkin diagram with only one non-zero label, by Lemma \ref{Slem1}, we may assume that $P$ is
either in $\odot^2\h$, or in $\Lambda^2\h$. The same is true for $\h=\sl(n,\Co)$ since it is given by the Dynkin diagram
with exactly two  non-zero labels, i.e.  $\h\otimes \h$ contains two irreducible $\h$-modules isomorphic to $\h$, one of
them coincides with $\P_1(\h)$ and we need to explore the other one.

If $P\in \Lambda^2\h$, then $P([x,y])=[P(x),y]+[x,P(y)]$, i.e. $P$ is a derivative of $\h$. Since $\h$ is simple, $P$ is
of the form $P(\cdot)=[\cdot,x]$ for some $x\in\h$, i.e. $P\in\P_1(\h)$.

Suppose that $P\in\odot^2 \h$, then $P([x,y])=-[P(x),y]-[x,P(y)]$. We have \begin{align*}
0=&P([[x,y],z]+[[y,z],x]+[[z,x],y])\\ =&-[P([x,y]),z]-[[x,y],P(z)]-[P([y,z]),x]-[[y,z],P(x)]-[P([z,x]),y]-[[z,x],P(y)]\\
=&[[P(x),y],z]+[[x,P(y)],z]-[[x,y],P(z)]+[[P(y),z],x]+[[y,P(z)],x]-[[y,z],P(x)]\\
&+[[P(z),x],y]+[[z,P(x)],y]-[[z,x],P(y)]\\ =&-2([[x,y],P(z)]+[[y,z],P(x)]+[[z,x],P(y)]).\end{align*} The last equality
implies $$[[P(x),y],z]+[[x,P(y)],z]+[[P(y),z],x]+[[y,P(z)],x]+[[P(z),x],y]+[[z,P(x)],y]=0.$$ Hence,
$[P([x,y]),z]+[P([y,z]),x]+[P([z,x]),y]=0$, i.e. $P([\cdot,\cdot])\in\R(\h)$. This shows that $P([x,y])=c[x,y]$ for all
$x,y\in\h$ and some constant $c\in\Co$. We conclude that $P=c\id_\h$. It is clear that $P\in\P(\h)$ if and only if $c=0$.
The lemma is proved.  \qed

\vskip0.3cm

We are left with the representations 1, 2, 4-8 from Table \ref{symBer}.

\vskip0.3cm

For the representations 4, 5, and 8 we prove that $\P_0(\h)=0$ using Lemma \ref{lemRP}. For each irreducible submodule
$V_\lambda\subset V\otimes \h$ different from the highest one and from $V$ we find a submodule $U\subset\Lambda^2
V\otimes\h$ such that the natural map $\tau$ from $(\Lambda^2 V\otimes\h)\otimes V$ to $V\otimes \h$, mapping $R\otimes x$
to $R(\cdot,x)$,  maps $U\otimes V$ onto $V_\lambda$. Since $U\not\subset\R(\h)$, we get that $V_\lambda\not\subset
\P(\h)$.

{\bf The subalgebra} $\h=\so(9,\Co)\subset\so((\Delta_9)^\Co)=\so(16,\Co)$. It holds $V\otimes\h=V\oplus
V_{\pi_2+\pi_4}\oplus V_{\pi_1+\pi_4}$. The submodule $ V_{\pi_2+\pi_4}\subset V\otimes \h$ is the highest one and we are
left with the submodule $V_{\pi_1+\pi_4}$. The $\h$-module $\odot^2\h\subset\Lambda^2 V\otimes\h$ contains the submodule
$V_{2\pi_1}$. We have $V_{2\pi_1}\otimes V=V_{2\pi_1+\pi_4}\oplus V_{\pi_1+\pi_4}$. Since $V_{2\pi_1+\pi_4}\not\subset
V\otimes\h$ and $\tau(V_{2\pi_1}\otimes V)\neq 0$, we get that $\tau(V_{2\pi_1}\otimes V)= V_{\pi_1+\pi_4}$. Since
$V_{2\pi_1}\not\subset\R(\h)$, we conclude that $V_{\pi_1+\pi_4}\not\subset\P(\h)$.

{\bf The subalgebra} $\h=\so(16,\Co)\subset\so((\Delta_{16}^+)^\Co)=\so(128,\Co)$. It holds $V\otimes\h=V\oplus
V_{\pi_2+\pi_8}\oplus V_{\pi_1+\pi_7}$.
 The submodule
$V_{\pi_2+\pi_8}\subset V\otimes \h$ is the highest one. The $\h$-module $\odot^2\h\subset\Lambda^2 V\otimes\h$ contains
the submodule $V_{2\pi_1}$. We have $V_{2\pi_1}\otimes V=V_{\pi_1+\pi_7}\oplus V_{2\pi_1+\pi_8}$. We conclude that
$V_{\pi_1+\pi_7}\not\subset\P(\h)$.

{\bf The subalgebra} $\h=\sp(8,\Co)\subset\so(V_{\pi_4})=\so(42,\Co)$. It holds $V\otimes\h=V\oplus V_{2\pi_1+\pi_4}\oplus
V_{\pi_1+\pi_3}$. The submodule $V_{2\pi_1+\pi_4}\subset V\otimes \h$ is the highest one. The $\h$-module
$\odot^2\h\subset\Lambda^2 V\otimes\h$ contains the submodule $V_{\pi_2}$. We have $V_{\pi_2}\otimes
V=V_{\pi_2+\pi_4}\oplus V_{\pi_1+\pi_3}\oplus V_{\pi_2}$. We conclude that $V_{\pi_1+\pi_3}\not\subset\P(\h)$.

\vskip0.5cm

The above trick does not work with the representations 1, 2, 6, 7 from Table \ref{symBer} and we use the direct
computations.

{\bf The subalgebra} $\h=\sp(2n,\Co)\subset\so(V_{\pi_2})$. It holds $V\otimes\h=V\oplus V_{2\pi_1+\pi_2}\oplus
V_{\pi_1+\pi_3}\oplus V_{2\pi_1}$. The submodule $ V_{2\pi_1+\pi_2}\subset V\otimes \h$ is the highest one and we are left
with the submodules $V_{\pi_1+\pi_3}$ and $V_{2\pi_1}$. The $\h$-module $\Lambda^2 V\otimes\h$ contains the submodule
$V_{\pi_2+\pi_4}$, in the same time $V_{\pi_2+\pi_4}\otimes V$  contains the submodule $V_{\pi_1+\pi_3}$ and it contains
none of the submodules $V_{2\pi_1+\pi_2}$ and $V_{2\pi_1}$.  We conclude that $V_{\pi_1+\pi_3}\not\subset\P(\h)$.

Consider the submodule $V_{2\pi_1}$. Let $e_{1},...,e_{n},e_{-1},...,e_{-n}$ be the standard basis of $\Co^{2n}$( such
that $\omega(e_i,e_{-i})=1$). The highest vector of the module $V_{2\pi_1}$ equals to $$\varphi=\sum_{i=2}^n e_1\wedge
e_i\otimes (E_{1,i}-E_{n+i,n+1}),$$ where $E_{a,b}$ is the matrix with $1$ on the position $(a,b)$ and zeros on the other
positions. To find $\varphi$, we consider the vector subspace of $V\otimes \h$ of weight $2\pi_1$ and find a vector
(defined up to a constant) annihilated by the generators of $\h$ corresponding to the simple positive roots. Let
$x=e_{-1}\wedge e_{-2}$, $y=e_{2}\wedge e_{3}$ and $z=e_{-1}\wedge e_{-3}$. Then
$ (\varphi(x)y,z)+ (\varphi(y)z,x)+ (\varphi(z)x,y)=2.$ Thus, $\varphi\not\in\P(\h)$ and $V_{2\pi_1}\not\subset
\P(\h)$.

{\bf The subalgebra} $\h=\so(n,\Co)\subset\so(V_{2\pi_1})$. It holds $V\otimes\h=V\oplus V_{2\pi_1+\pi_2}\oplus
V_{\pi_1+\pi_3}\oplus V_{\pi_2}$. The submodule $ V_{2\pi_1+\pi_2}\subset V\otimes \h$ is the highest one and we are left
with the submodules $V_{\pi_1+\pi_3}$ and $V_{\pi_2}$. The $\h$-module $\odot^2\h\subset\Lambda^2 V\otimes\h$ contains
the submodule $V_{\pi_4}$, and $V_{\pi_4}\otimes V$ contains $V_{\pi_1+\pi_3}$ and it  contains none of the submodules
$V_{2\pi_1+\pi_2}$ and $V_{\pi_2}$.  We conclude that $V_{\pi_1+\pi_3}\not\subset\P(\h)$.

Consider the submodule $V_{\pi_2}$. Let $e_{1},...,e_{m},e_{-1},...,e_{-m}$ and $e_{1},...,e_{m},e_{-1},...,e_{-m},e_0$ be
the standard bases of $\Co^{2m}$ and $\Co^{2m+1}$, respectively ( such that $g(e_i,e_{-i})=1$ and $g(e_0,e_{0})=1$). The
highest vector of the module $V_{\pi_2}$ equals to
\begin{align*} \varphi=&\sum_{i=3}^m e_1\odot e_i\otimes (E_{2,i}-E_{m+i,m+2})-\sum_{i=3}^m e_2\odot e_i\otimes (E_{1,i}-E_{m+i,m+1}),\quad \text { if } n=2m,\\
\varphi=&e_1\odot e_0\otimes (E_{2,2m+1}-E_{2m+1,m+2})-e_2\odot e_0\otimes (E_{1,2m+1}-E_{2m+1,m+1})\\&+ \sum_{i=3}^m
e_1\odot e_i\otimes (E_{2,i}-E_{m+i,m+2})-\sum_{i=3}^m e_2\odot e_i\otimes (E_{1,i}-E_{m+i,m+1}),\quad  \text { if }
n=2m+1.\end{align*} Taking in the both cases  $x=e_{-1}\odot e_{-3}$, $y=e_{1}\odot e_{3}$ and $z=e_{-1}\odot e_{-2}$, we
get $ (\varphi(x)y,z)+ (\varphi(y)z,x)+ (\varphi(z)x,y)=1.$ Thus, $\varphi\not\in\P(\h)$ and
$V_{\pi_2}\not\subset \P(\h)$.

{\bf The subalgebra} $\sl(8,\Co)\subset \so(\Lambda^4\Co^{8})$. We have $V\otimes\g=V\oplus V_{\pi_1+\pi_4+\pi_7}\oplus
V_{\pi_1+\pi_3}\oplus V_{\pi_5+\pi_7}$. The submodule $ V_{\pi_1+\pi_4+\pi_7}\subset V\otimes \h$ is the highest one and
we are left with the submodules $V_{\pi_1+\pi_3}$ and $V_{\pi_5+\pi_7}$. Let $e_1,...,e_8$ be the standard basis of
$\Co^8$. The metric on $V=\Lambda^4\Co^8$ is given by the exterior multiplication,
$ (\omega,\theta)=\omega\wedge\theta=\theta\wedge\omega\in \Lambda^8\Co^8\simeq \Co$, we assume that $ (e_1\wedge
e_2\wedge e_3\wedge e_4,e_5\wedge e_6\wedge e_7\wedge e_8)=1$. The highest vector of the submodule $V_{\pi_1+\pi_3}\subset
V\otimes\h$ equals to $\varphi=\sum_{i=1}^5 e_1\wedge e_2\wedge e_3\wedge e_i\otimes E_{1,i}$. Taking $x=e_5\wedge
e_6\wedge e_7\wedge e_8$, $y=e_2\wedge e_4\wedge e_5\wedge e_6$, and $z=e_3\wedge e_4\wedge e_7\wedge e_8$, we get
$ (\varphi(x)y,z)+ (\varphi(y)z,x)+ (\varphi(z)x,y)=-1.$ Hence, $V_{\pi_1+\pi_3}\not\subset \P(\h)$. The symmetry
of the Dynkin diagram of $\sl(8,\Co)$ implies $V_{\pi_5+\pi_7}\not\subset \P(\h)$.

{\bf The subalgebra} $\h=F_4^\Co\subset\so(26,\Co)$. To deal with this representation we use the following description of
it form \cite{Adams}. The Lie algebra $F^\Co_4$ admits the structure of $\mathbb{Z}_2$-graded Lie algebra:
$F^\Co_4=\so(9,\Co)\oplus(\Delta_9)^\Co$. The representation space $\Co^{26}$ is decomposed into the direct sum
$\Co^{26}=\Co\oplus\Co^9\oplus(\Delta_9)^\Co$. The elements of the subalgebra $\so(9,\Co)\subset F^\Co_4$ preserve these
components, annihilate $\Co$ and act naturally on $\Co^9$ and $(\Delta_9)^\Co$. Elements of $(\Delta_9)^\Co\subset
F^\Co_4$ take $\Co$ and $\Co^9$ to $(\Delta_9)^\Co$ (multiplication by constants and the Clifford multiplication,
respectively), and take $(\Delta_9)^\Co$ to $\Co\oplus\Co^9$ (the charge conjugation plus the natural map assigning a
vector to a pair of spinors). Let $P\in\P(F^\Co_4)$. Decompose it as the sum $P=\varphi+\psi$, where $\varphi$ and $\psi$
take values in $\so(9,\Co)$ and $(\Delta_9)^\Co$, respectively. The condition $P\in\P(F^\Co_4)$ implies
\begin{align*}
&\varphi|_{\Co^9}\in\P(\so(9,\Co)\subset\so(9,\Co)),\quad
\varphi|_{(\Delta_9)^\Co}\in\P(\so(9,\Co)\subset\so((\Delta_9)^\Co)),
\\& \varphi(a)=0,\quad (\psi(a)x,s)+(\psi(x)s,a)=0\\ &(\psi(x)y-\psi(y)x,s)+(\varphi(s)x,y)=0,\quad
(\psi(r)s-\psi(s) r,x)+(\varphi(x)r,s)=0.
\end{align*} for all $a\in\Co$, $x,y\in\Co^9$, and $s,r\in (\Delta_9)^\Co$.

We will denote the $\h$-modules by $V^{\h}_\lambda$ and the $\so(9,\Co)$-modules by $V_{\lambda}$. We have $V\otimes
\h=V\oplus V^{\h}_{\pi_1+\pi_4}\oplus V^{\h}_{\pi_2}$. The submodule $V^{\h}_{\pi_1+\pi_4}\subset V\otimes\h$ is the
highest one and we need to explore the module $V^{\h}_{\pi_2}$. Note that $\dim V^{\h}_{\pi_2}=273$. The above equalities
show that $P$ is uniquely defined by $\psi|_{\Co^9\oplus (\Delta_9)^\Co}.$ In particular, $$\psi|_{(\Delta_9)^\Co}\in
(\Delta_9)^\Co\otimes (\Delta_9)^\Co=\Co\oplus V_{2\pi_4}\oplus V_{\pi_3}\oplus V_{\pi_2}\oplus V_{\pi_1}$$ defines
$\varphi|_{\Co^9}\in\P(\so(9,\Co)\subset\so(9,\Co))=V_{\pi_1}\oplus V_{\pi_1+\pi_2}$. Hence, $\varphi|_{\Co^9}\in
V_{\pi_1}=\P_1(\so(9,\Co)\subset\so(9,\Co))$. Next, $$\psi|_{\Co^9}\in \Co^9\otimes (\Delta_9)^\Co=V_{\pi_4}\oplus
V_{\pi_1+\pi_4}$$ defines $\varphi|_{(\Delta_9)^\Co}\in\P(\so(9,\Co)\subset\so((\Delta_9)^\Co)=V_{\pi_4}$.

It is clear that $\P_1(\h)$ is given by $\Co\oplus V_{\pi_1}\subset (\Delta_9)^\Co\otimes (\Delta_9)^\Co$ and by
$V_{\pi_4}\subset \Co^9\otimes (\Delta_9)^\Co$ (recall that $\Co\oplus V_{\pi_1}\oplus V_{\pi_4}=\Co\oplus\Co^9\oplus
(\Delta_9)^\Co=V$). The dimensions of the $\so(9,\Co)$-modules $V_{2\pi_4}$, $V_{\pi_3}$ $V_{\pi_2}$  $V_{\pi_1}$, and
$V_{\pi_1+\pi_4}$ equal, respectively, 126, 84, 36, 9, and 128. We see that the sum of some of these numbers can not equal
$\dim V^{\h}_{\pi_2}=273$. This shows that $V^{\h}_{\pi_2}\not\subset\P(\h)$.

{\bf Acknowledgment.} I am thankful to D.\,V.\,Alekseevsky and Thomas Leistner for useful and stimulating discussions on the topic of this paper. The work was supported by the grant  201/09/P039 of the
Grant Agency of Czech Republic   and by the grant MSM~0021622409 of the Czech Ministry of Education.

\end{document}